\documentclass[final,3p,times]{elsarticle}

\usepackage{amssymb,amsmath,amsfonts,amsthm,mathrsfs}
\usepackage{diagbox}

\newdefinition{dfn}{Definition}
\newdefinition{ex}{Example}
\newtheorem{thm}{Theorem}
\newtheorem{lmm}{Lemma}

\journal{Statistics \& Probability Letters}

\makeatletter \def\ps@pprintTitle{ \let\@oddhead\@empty \let\@evenhead\@empty \def\@oddfoot{\reset@font\hfil} \let\@evenfoot\@oddfoot} \makeatother

\begin{document}

\begin{frontmatter}

\title{Universal distribution of the empirical coverage in split conformal prediction}

\author{Paulo C. Marques F.}

\ead{PauloCMF1@insper.edu.br}

\affiliation{
  organization = {Insper Institute of Education and Research},
  addressline = {Rua Quatá, 300},
  city = {São Paulo},
  postcode = {04546-042}, 
  state = {SP},
  country = {Brazil}
}

\begin{abstract}
When split conformal prediction operates in batch mode with exchangeable data, we determine the exact distribution of the empirical coverage of prediction sets produced for a finite batch of future observables, as well as the exact distribution of its almost sure limit when the batch size goes to infinity. Both distributions are universal, being determined solely by the nominal miscoverage level and the calibration sample size, thereby establishing a criterion for choosing the minimum required calibration sample size in applications.
\end{abstract}



\end{frontmatter}


\section{Introduction}\label{sec:intro}

Conformal prediction is a framework developed to quantify the confidence in the forecasts made by general predictive models which is quickly moving the field of machine learning from a stage dominated by point predictions to a new period in which forecasts are summarized by prediction sets with statistical guarantees. Several features make conformal prediction appealing for use with contemporary machine learning algorithms: it is universal (distribution-free), able to handle high-dimensional data, model agnostic, and its properties hold for finite samples \cite{vovk1999,vovk2005,fontana}. Notably, the split conformal prediction algorithm \cite{papadopoulos,lei} is a widely adopted conformalization technique which strikes a balance between predictive properties and computational cost. Our goal in this paper is to identify the exact distribution of the empirical coverage of prediction sets produced by the split conformal prediction procedure for a finite batch of future observables, as well as to determine the exact distribution of its almost sure limit when the batch size tends to infinity. Both distributions are universal in the sense that they are determined solely by the nominal miscoverage level and the calibration sample size. The distribution of the empirical coverage was investigated for the first time in \cite{vovk2013} and \cite{hulsman}, with further discussion in \cite{angelopoulos}. Our contribution consists in a formulation that emphasizes the role of the data exchangeability assumption and the combinatorial nature of the aforementioned properties of the empirical coverage, which are derived using standard exchangeability tools. Although this investigation pertains to the foundations of conformal prediction, the results are eminently practical and lead to a criterion summarized in Table 1 for the choice of the calibration sample size in applications.

\section{Split conformal prediction}\label{sec:splitcp}

Let $(\Omega,\mathscr{F},P)$ denote the underlying probability space from which we induce the distributions of all random objects considered in the paper.

\begin{dfn}
A sequence of random objects $\{O_i\}_{i\geq 1}$ is {\it exchangeable} if, for every $n\geq 1$ and every permutation $\pi~:~\{1,\dots,n\} \stackrel{\cong}{\longrightarrow} \{1,\dots,n\}$, the random tuples $(O_1,\dots,O_n)$ and $(O_{\pi(1)},\dots,O_{\pi(n)})$ have the same distribution.
\end{dfn}

We are in a supervised learning setting \cite{esl} in which for each sample unit we have a $d$-dimensional vector of predictors $X_i\in\mathbb{R}^d$ and a response variable $Y_i\in\mathscr{Y}$. Specifically, in regression tasks with univariate response we take $\mathscr{Y}=\mathbb{R}$ and in classification problems $\mathscr{Y}=\{1,\dots,L\}$ is a set of class labels. We have a data sequence of random pairs laid out as
\begin{equation*}\label{eq:seq}
  \underbrace{(X_{-t+1},Y_{-t+1}),\dots,(X_0,Y_0)}_\text{training},\underbrace{(X_1,Y_1),\dots,(X_n,Y_n)}_\text{calibration},\underbrace{(X_{n+1},Y_{n+1}),\dots}_\text{future}
\end{equation*}
which is modeled by us as being exchangeable. At the beginning of the sequence we have the {\it training sample} $T=((X_{-t+1},Y_{-t+1}),\dots,(X_0,Y_0))$, of size $t\geq 1$, followed by the {\it calibration sample} $((X_1,Y_1),\dots,(X_n,Y_n))$, of size $n\geq 1$, and the sequence of {\it future observables} $\{(X_{n+i},Y_{n+i})\}_{i\geq 1}$. In applications, the available data is randomly split into the training and calibration samples, hence the name {\it split} conformal prediction, also known as the inductive case of conformal prediction. The data exchangeability assumption allows us to conveniently place the training sample at the beginning of the sequence. Let $\mathscr{T} = \sigma(T)$ be the smallest sub-$\sigma$-field of $\mathscr{F}$ with respect to which the training sample is measurable. 

\begin{dfn}\label{dfn:score}
A {\it conformity function} is a mapping $\rho:\mathbb{R}^d\times\mathscr{Y}\times\Omega\to\mathbb{R}$ such that $\rho(x,y)=\rho(x,y,\,\boldsymbol{\cdot}\;)$ is $\mathscr{T}$-measurable for every $x\in\mathbb{R}^d$ and every $y\in\mathscr{Y}$. The sequence of {\it conformity scores} $\{S_i\}_{i\geq 1}$ associated with a conformity function $\rho$ is defined by $S_i(\omega)=\rho(X_i(\omega),Y_i(\omega),\omega)$. We say that a conformity function $\rho$ is {\it regular} with respect to a specific data sequence if there are no ties among the corresponding conformity scores $\{S_i\}_{i\geq 1}$ almost surely.
\end{dfn}

Note that the regularity of a specific conformity function $\rho$ is contextual, being inherently dependent on the distribution of the underlying data sequence. Technically, we can always avoid ties among the sequence of conformity scores almost surely by introducing a properly constructed ancillary tie-breaking sequence.

\begin{ex}\label{ex:scores}
In regression problems, let $\hat{\mu}:\mathbb{R}^d\times\Omega\to\mathbb{R}$ be a regression function estimator. A standard choice \cite{papadopoulos,lei} is to use the conformity function $\rho(x,y,\omega)=|y-\hat{\mu}(x,\omega)|$. Conformalized quantile regression \cite{romano} is a widely adopted alternative. For $p\in[0,1]$, let $\xi_p(x)=\inf\,\{y\in\mathbb{R}:P(Y_1\leq y\mid X_1=x\}\geq p\}$ be the conditional $p$th quantile function and suppose that we have an estimator $\hat{\xi}_p:\mathbb{R}^d\times\Omega\to\mathbb{R}$ of $\xi_p$. Choose $0<p_\text{lo}\leq p_\text{hi}<1$ and define the conformity function $\rho(x,y,\omega) = \max\,\{\hat{\xi}_{p_\text{hi}}(x,\omega)-y,y-\hat{\xi}_{p_\text{lo}}(x,\omega)\}$. The choice of $p_\text{lo}$ and $p_\text{hi}$ is discussed in \cite{romano}. In classification problems, let $\hat{\pi}_\ell:\mathbb{R}^d\times\Omega\to[0,1]$ be a classification algorithm outputting probabilities for each one of the class labels $\ell=1,\dots,L$. In this classification case, we can take our conformity function to be $\rho(x,y,\omega) = 1-\hat{\pi}_{y}(x,\omega)$.
\end{ex}

Conformity functions are agnostic to the choice of the specific models or algorithms used to construct $\hat{\mu}$, $\hat{\xi}_p$, and $\hat{\pi}$ in Example \ref{ex:scores}. The intuition is that the associated conformity scores measure the ability of the model to make accurate predictions on the calibration sample, whose information is not used in the model´s training process, and the assumed data sequence exchangeability transfers this assessment of the model's predictive capacity from the calibration sample to the sequence of future observables. The following result is proved in the Appendix.

\begin{lmm}\label{lmm:xchscr}
Under the data exchangeability assumption, the sequence of conformity scores $\{S_i\}_{i\geq 1}$ is exchangeable.
\end{lmm}

For a real number $t$, let $\lceil t\rceil=\min\,\{k\in\mathbb{Z}:t\leq k\}$ and $\lfloor t\rfloor=\max\,\{k\in\mathbb{Z}:k\leq t\}$ denote the ceiling and the floor of $t$, respectively. 

\begin{dfn}\label{dfn:cpset}
For a regular conformity function $\rho$, denote the associated ordered calibration sample conformity scores by $S_{(1)},S_{(2)},\dots,S_{(n)}$. Let $0<\alpha<1$ be a specified {\it nominal miscoverage level} satisfying $\lceil(1-\alpha)(n+1)\rceil\leq n$, in which case we say that the pair $(n,\alpha)$ is {\it feasible}. Define the random {\it conformal prediction set} $D^{(\alpha)}_n:\mathbb{R}^d\times\Omega\to\mathscr{A}$ by
\begin{equation*}
  D^{(\alpha)}_n(x,\omega) = \left\{y\in\mathscr{Y}:\rho(x,y,\omega) \leq S_{(\lceil(1-\alpha)(n+1)\rceil)}(\omega)\right\},
\end{equation*}
for a suitable $\sigma$-field of subsets of $\mathscr{Y}$ denoted by $\mathscr{A}$. We use the notation $D^{(\alpha)}_n(x)=D^{(\alpha)}_n(x,\,\boldsymbol{\cdot}\;)$.
\end{dfn}

\begin{ex}\label{ex:intervals}
Let $\hat{s}=S_{(\lceil(1-\alpha)(n+1)\rceil)}(\omega_0)$, for an outcome $\omega_0\in\Omega$. It follows from Definition \ref{dfn:cpset} and the definitions in Example \ref{ex:scores} that for a future vector of predictors $x^*\in\mathbb{R}^d$ the observed conformal prediction sets have the forms: $D^{(\alpha)}_n(x^*)~=~(\hat{\mu}(x^*)-\hat{s},\hat{\mu}(x^*)+\hat{s})$, for the standard conformity function, $D^{(\alpha)}_n(x^*)=(\hat{\xi}_{p_\text{lo}}(x^*)-\hat{s},\hat{\xi}_{p_\text{hi}}(x^*)+\hat{s})$, for conformalized quantile regression, and $D^{(\alpha)}_n(x^*)=\{y\in\{1,\dots,L\}:\hat{\pi}_y(x^*)\geq 1-\hat{s}\}$, for classification.
\end{ex}

The first major consequence of Lemma \ref{lmm:xchscr} is the classical marginal validity property \cite{papadopoulos,lei} of the conformal prediction sets introduced in Definition \ref{dfn:cpset}. This property can be described briefly as follows. For a regular conformity function $\rho$, using the notations and conditions in Definition \ref{dfn:cpset}, the distributional symmetry expressed in Lemma~\ref{lmm:xchscr} and the fact that we have no ties almost surely among the sequence of conformity scores, ensure that, for every $i\geq 1$, the conformity score $S_{n+i}$ for a future random pair $(X_{n+i},Y_{n+i})$ is uniformly ranked among the ordered calibration sample conformity scores: $P(S_{n+i}\leq S_{(j)})=j/(n+1)$, for $j=1,\dots,n$. By choosing $j=\lceil(1-\alpha)(n+1)\rceil$, noting that $t\leq\lceil t\rceil<t+1$, for every $t\in\mathbb{R}$, and considering that $Y_{n+i} \in D^{(\alpha)}_n(X_{n+i})$ if and only if $S_{n+i}\leq S_{(\lceil(1-\alpha)(n+1)\rceil)}$ (as per Definitions \ref{dfn:score} and \ref{dfn:cpset}), we have the marginal validity property:
\begin{equation}\label{eq:mvp}
  1 - \alpha \leq P(Y_{n+i} \in D^{(\alpha)}_n(X_{n+i})) < 1 - \alpha + \frac{1}{n+1}\,.
\end{equation}

\section{Empirical coverage distribution}\label{sec:ecd}

\begin{dfn}\label{dfn:coverage}
Using the notations and conditions in Definition \ref{dfn:cpset}, let $\{Z_i\}_{i\geq 1}$ be a sequence of {\it coverage indicators} defined by $Z_i=1$, if $Y_{n+i}\in D^{(\alpha)}_n(X_{n+i})$, and $Z_i=0$, otherwise. The {\it empirical coverage} of a batch of $m\geq 1$ future observables is the random variable $C^{(n,\alpha)}_m=(1/m)\sum_{i=1}^m Z_i$.
\end{dfn}

In general, the coverage indicators $Z_i$ are dependent random variables, since for all future observables the corresponding conformal prediction sets in Definition \ref{dfn:cpset} are defined in terms of the same calibration sample conformity score $S_{(\lceil(1-\alpha)(n+1)\rceil)}$. This would still be the case even if we had started with the stronger assumption of an independent and identically distributed data sequence. The interesting fact is that Definition \ref{dfn:coverage} inherits through Lemma \ref{lmm:xchscr} the distributional symmetry implied by the data exchangeability assumption, giving us the following result, proved in the Appendix.

\begin{thm}\label{thm:finite}
Under the data exchangeability assumption, for a regular conformity function, the sequence of coverage indicators $\{Z_i\}_{i\geq 1}$ is exchangeable and $m\times C^{(n,\alpha)}_m$ is distributed as a $\text{Beta-Binomial}(\lceil(1-\alpha)(n+1)\rceil,\lfloor\alpha(n+1)\rfloor)$ random variable, to the effect that the distribution of the empirical coverage is given by
$$
  P\left(C^{(n,\alpha)}_m=\frac{k}{m}\right) = \binom{m}{k}\,\frac{n!\,(k+\lceil(1-\alpha)(n+1)\rceil-1)!\,(m-k+\lfloor\alpha(n+1)\rfloor-1)!}{(\lceil(1-\alpha)(n+1)\rceil-1)!\,(\lfloor\alpha(n+1)\rfloor-1)!\,(m+n)!}\,,
$$
for $k=0,1,\dots,m$, and every future batch size $m\geq 1$.
\end{thm}

By symmetry, Definition \ref{dfn:coverage} and the exchangeability of the sequence of coverage indicators established in Theorem \ref{thm:finite} yield that $\textrm{E}[C^{(n,\alpha)}_m]=\textrm{E}[Z_1]=P(Y_{n+1} \in D^{(\alpha)}_n(X_{n+1}))$. Consequently, we can interpret the marginal validity property (\ref{eq:mvp}) as partial information about the distribution of the empirical coverage. Specifically, as an inequality constraint on the expectation of $C^{(n,\alpha)}_m$. The following result, proved in the Appendix as a direct consequence of Theorem \ref{thm:finite} and de Finetti's representation theorem, identifies the distribution of the almost sure limit of the empirical coverage $C^{(n,\alpha)}_m$ when the future batch size $m$ tends to infinity.

\begin{thm}\label{thm:limit}
Under the data exchangeability assumption, for a regular conformity function, the empirical coverage $C^{(n,\alpha)}_m$ converges almost surely, when the future batch size tends to infinity, to a random variable $C^{(n,\alpha)}_\infty$ with distribution $\text{Beta}(\lceil(1-\alpha)(n+1)\rceil,\lfloor\alpha(n+1)\rfloor)$.
\end{thm}

\section{Calibration sample size and concluding remarks}\label{sec:calsize}

In applications, we typically use a trained model to construct a large number of prediction intervals and Theorem~\ref{thm:limit} gives us a criterion to determine the minimum calibration sample size required to control the empirical coverage $C^{(n,\alpha)}_\infty$ of an infinite batch of future observables. Given a nominal miscoverage level $\alpha$, we specify an $\epsilon>0$ and a tolerance probability $0<\tau<1$, looking for the smallest calibration sample size such that the empirical coverage $C^{(n,\alpha)}_\infty$ of an infinite batch of future observables is within $\epsilon$ of $1-\alpha$, with a probability of at least $\tau$. Formally, recalling from Definition~\ref{dfn:cpset} that a feasible pair $(n,\alpha)$ satisfies $\lceil(1-\alpha)(n+1)\rceil\leq n$, which is equivalent to saying that the integer $n\geq (1-\alpha)/\alpha$, the minimum required calibration size is given by
\begin{equation*}
  n_0 = \min\left\{n\geq (1-\alpha)/\alpha : P\left(\left|C^{(n,\alpha)}_\infty - (1 - \alpha)\right| < \epsilon\right) \geq \tau\right\}.
\end{equation*}
Table \ref{tab:uct} gives the values of the minimum required calibration sample size $n_0$ for different values of $\alpha$, $\epsilon$, and $\tau$. That an understanding of the distribution of the empirical coverage is necessary to determine the minimum required calibration sample sizes in applications of split conformal prediction was first discussed in \cite{angelopoulos}. The calibration sample sizes presented in \cite{angelopoulos} are slightly larger than the corresponding values in Table \ref{tab:uct}. In the repository \cite{githubcoverage} we have the \texttt{R} \cite{R} code used to determine the calibration sample sizes in Table \ref{tab:uct} and a comparison with the corresponding values given in \cite{angelopoulos}. Repository \cite{githubcoverage} also contains a simulation illustrating the results in Theorems \ref{thm:finite} and \ref{thm:limit}.

\begin{table}[t!]
\centering
\small
{\renewcommand{\arraystretch}{1.25}
\begin{tabular}{c|ccc|ccc|ccc|ccc|}
\cline{2-13}
                           & \multicolumn{3}{c|}{$\epsilon=0.1$} & \multicolumn{3}{c|}{$\epsilon=0.05$} & \multicolumn{3}{c|}{$\epsilon=0.01$} & \multicolumn{3}{c|}{$\epsilon=0.005$} \\ \hline
\multicolumn{1}{|c|}{\backslashbox{$1-\alpha$}{$\tau$}}  & 90\%       & 95\%       & 99\%      & 90\%       & 95\%       & 99\%       & 90\%       & 95\%       & 99\%       & 90\%        & 95\%       & 99\%       \\ \hline
\multicolumn{1}{|c|}{80\%} & 40         & 57         & 98        & 170        & 241        & 418        & 4,326      & 6,142      & 10,611     & 17,314      & 24,581     & 42,457     \\
\multicolumn{1}{|c|}{85\%} & 30         & 42         & 77        & 134        & 189        & 330        & 3,446      & 4,893      & 8,451      & 13,794      & 19,587     & 33,830      \\
\multicolumn{1}{|c|}{90\%} & 11         & 14         & 47        & 90         & 128        & 227        & 2,429      & 3,448      & 5,958      & 9,733       & 13,821     & 23,875     \\
\multicolumn{1}{|c|}{95\%} & 19         & 19         & 29        & 22         & 29         & 97         & 1,270      & 1,806      & 3,132      & 5,125       & 7,278      & 12,578     \\ \hline
\end{tabular}}
\caption{Universal coverage tolerance table for split conformal prediction. For a nominal miscoverage level $\alpha$, an $\epsilon>0$, and a tolerance probability $\tau$, the table entries are the minimum required calibration sample sizes such that the empirical coverage $C^{(n,\alpha)}_\infty$ of an infinite batch of future observables is within $\epsilon$ of $1-\alpha$, with a probability of at least $\tau$, according to Theorem \ref{thm:limit}.}\label{tab:uct}
\end{table}

\section*{Acknowledgments}

Paulo C. Marques F. receives support from FAPESP (Fundação de Amparo à Pesquisa do Estado de São Paulo) through project 2023/02538-0.

\appendix

\section*{Appendix. Proofs}\label{app:proofs}

\begin{proof}[Proof of Lemma \ref{lmm:xchscr}]
Recall that $T(\omega) = ((X_{-t+1}(\omega),Y_{-t+1}(\omega)),\dots,(X_0(\omega),Y_0(\omega)))$ and $\mathscr{T}=\sigma(T)$. Since the conformity function $\rho$ in Definition \ref{dfn:score} is such that $\rho(x,y)$ is $\mathscr{T}$-measurable for every $x\in\mathbb{R}^d$ and every $y\in\mathscr{Y}$, Doob-Dynkin's lemma (see \cite{schervish}, Theorem A.42) implies that there is a measurable function
$$
  h : \underbrace{(\mathbb{R}^d\times\mathscr{Y}) \times \dots \times (\mathbb{R}^d\times\mathscr{Y})}_{t+1 \text{ times}} \to \mathbb{R}
$$
such that $\rho(X_i(\omega),Y_i(\omega),\omega)=h(T(\omega),(X_i(\omega),Y_i(\omega)))$, for every $i\geq 1$, and each $\omega\in\Omega$. Hence, for Borel sets $B_1,\dots,B_n$, we have
\begin{align*}
  P\left(\cap_{i=1}^n \big\{\omega:S_i(\omega)\in B_i\big\}\right) &= P\left(\cap_{i=1}^n \big\{\omega:\rho(X_i(\omega),Y_i(\omega),\omega)\in B_i\big\}\right) \\
  &= P\left(\cap_{i=1}^n \big\{\omega:h(T(\omega),(X_i(\omega),Y_i(\omega)))\in B_i\big\}\right) \\
  &=P\left(\cap_{i=1}^n \big\{\omega:(T(\omega),(X_i(\omega),Y_i(\omega)))\in h^{-1}(B_i)\big\}\right). \tag{$\star$}
\end{align*}
For any permutation $\pi:\{-t+1,\dots,0,1,\dots,n\} \stackrel{\cong}{\longrightarrow} \{-t+1,\dots,0,1,\dots,n\}$,
define
$$
  T_\pi(\omega) = ((X_{\pi(-t+1)}(\omega),Y_{\pi(-t+1)}(\omega)),\dots,(X_{\pi(0)}(\omega),Y_{\pi(0)}(\omega))).
$$
If we consider only permutations $\pi$ such that $\pi(j)=j$, for $-t+1\leq j\leq 0$, then $T_\pi=T$ and the data exchangeability assumption yields
\begin{align*}
  (\star) &= P\left(\cap_{i=1}^n \big\{\omega:(T_\pi(\omega),(X_{\pi(i)}(\omega),Y_{\pi(i)}(\omega)))\in h^{-1}(B_i)\big\}\right) \\
  &= P\left(\cap_{i=1}^n \big\{\omega:(T(\omega),(X_{\pi(i)}(\omega),Y_{\pi(i)}(\omega)))\in h^{-1}(B_i)\big\}\right) \\
  &=P\left(\cap_{i=1}^n \big\{\omega:h(T(\omega),(X_{\pi(i)}(\omega),Y_{\pi(i)}(\omega)))\in B_i\big\}\right) \\
  &= P\left(\cap_{i=1}^n \big\{\omega:\rho(X_{\pi(i)}(\omega),Y_{\pi(i)}(\omega),\omega)\in B_i\big\}\right) \\
  &=P\left(\cap_{i=1}^n \big\{\omega:S_{\pi(i)}(\omega)\in B_i\big\}\right).
\end{align*}
Since this restriction on $\pi$ still allows an arbitrary permutation of the conformity scores $S_{1},\dots,S_{n}$, and the argument holds for every $n\geq 1$, the desired exchangeability of the sequence of conformity scores $\{S_i\}_{i\geq 1}$ follows.
\end{proof}

\begin{proof}[Proof of Theorem \ref{thm:finite}]
Let $b=\lceil(1-\alpha)(n+1)\rceil$ and $g=n-b+1=\lfloor\alpha(n+1)\rfloor$, recalling from Definition \ref{dfn:cpset} that the pair $(n,\alpha)$ is assumed to be feasible, so that $b\leq n$. We will prove by induction on the batch size $m\geq 1$ that
\begin{equation*}
  P(Z_1=z_1,\dots,Z_m=z_m) = \frac{n!\,(k+b-1)!\,(m-k+g-1)!}{(b-1)!\,(g-1)!\,(m+n)!}\,, \tag{$\ast$}
\end{equation*}
in which $(z_1,\dots,z_m)\in\{0,1\}^m$, with $\sum_{i=1}^m z_i=k$, for $k=0,1,\dots,m$. Due to the assumed regularity of the underlying conformity function $\rho$, there are no ties among the sequence of conformity scores almost surely, and the distributional symmetry established in Lemma \ref{lmm:xchscr} implies that $S_{n+1}$ is uniformly ranked among the calibration sample conformity scores $(S_{1},\dots,S_{n})$. By Definitions \ref{dfn:cpset} and \ref{dfn:coverage}, for every $i\geq 1$, the coverage indicator $Z_i=1$ if and only if $S_{n+i}\leq S_{(b)}$, so that $P(Z_1=1)=P(S_{n+1}\leq S_{(b)})=b/(n+1)$. Hence, $P(Z_1=0)=g/(n+1)$ and property $(\ast)$ holds for $m=1$. By Lemma \ref{lmm:xchscr} and the regularity of $\rho$, the conformity score $S_{n+m+1}$ is ranked uniformly among the $n+m$ conformity scores $(S_1,\dots,S_n,S_{n+1},\dots,S_{n+m})$. Moreover, the event $\{Z_1=z_1,\dots,Z_m=z_m\}$, with $\sum_{i=1}^m z_i=k$, means that exactly $k$ of the conformity scores $(S_{n+1},\dots,S_{n+m})$ are less than or equal to $S_{(b)}$.  Hence, given that $\{Z_1=z_1,\dots,Z_m=z_m\}$, with $\sum_{i=1}^m z_i=k$, we have that $S_{n+m+1}\leq S_{(b)}$ if and only if $S_{n+m+1}$ is ranked among the $b+k$ conformity scores $\{S_{(1)},\dots,S_{(b)}\} \cup \left(\cup_{i=1}^m \left\{S_{n+i} : S_{n+i}\leq S_{(b)}\right\}\right)$, to the effect that
\begin{align*}
  P(Z_{m+1}=1\mid Z_1=z_1,\dots,Z_m=z_m) &= P(S_{n+m+1}\leq S_{(b)}\mid Z_1=z_1,\dots,Z_m=z_m) \\
  &= \frac{b+k}{m+n+1}\,.
\end{align*}
Now, for the inductive step, suppose that property $(*)$ holds for some batch size $m\geq 2$. The product rule gives
\begin{align*}
  P(Z_1=z_1,\dots,Z_m=z_m,Z_{m+1}=1) &= P(Z_{m+1}=1\mid Z_1=z_1,\dots,Z_m=z_m) \times P(Z_1=z_1,\dots,Z_m=z_m) \\ &= \frac{n!\,(k+b)!\,(m-k+g-1)!}{(b-1)!\,(g-1)!\,(m+n+1)!}\,.
\end{align*}
Since $P(Z_1=z_1,\dots,Z_m=z_m,Z_{m+1}=0)=1-P(Z_1=z_1,\dots,Z_m=z_m,Z_{m+1}=1)$, in general we have that
$$
  P(Z_1=z_1,\dots,Z_{m+1}=z_{m+1}) = \frac{n!\,(k'+b-1)!\,(m-k'+g)!}{(b-1)!\,(g-1)!\,(m+n+1)!}
$$
in which $(z_1,\dots,z_{m+1})\in\{0,1\}^{m+1}$, with $\sum_{i=1}^{m+1} z_i=k'$, for $k'=0,1,\dots,m+1$. Therefore, property $(\ast)$ holds for a batch with size $m+1$, completing the inductive step and implying that property $(\ast)$ holds for every batch size $m\geq 1$. Inspection of the right hand side of $(\ast)$ reveals that the random vector $(Z_1,\dots,Z_m$) is exchangeable, and since this holds for every batch size $m\geq 1$, we get as our first conclusion that the sequence of coverage indicators $\{Z_i\}_{i\geq 1}$ is exchangeable. Finally, the event $\{\sum_{i=1}^m Z_i=k\}$ is the union of $\binom{m}{k}$ mutually exclusive and, by exchangeability, equiprobable events of the form $\{Z_1=z_1,\dots,Z_m=z_m\}$, in which $(z_1,\dots,z_m)\in\{0,1\}^m$, with $\sum_{i=1}^m z_i=k$. Therefore, property $(\ast)$ and Definition \ref{dfn:coverage} yield the desired result:
$$
  P\left(C^{(n,\alpha)}_m=\frac{k}{m}\right) = P\!\left({\textstyle \sum_{i=1}^m Z_i=k}\right) = \binom{m}{k}\,\frac{n!\,(k+b-1)!\,(m-k+g-1)!}{(b-1)!\,(g-1)!\,(m+n)!}\,.
$$
\end{proof}

\begin{proof}[Proof of Theorem \ref{thm:limit}]
By Theorem \ref{thm:finite}, the sequence of coverage indicators $\{Z_i\}_{i\geq 1}$ is exchangeable, and de Finetti's representation theorem \cite{schervish} states that there is a random variable, say, $C^{(n,\alpha)}_\infty:\Omega\to[0,1]$, with distribution $\mu$, such that, given that $C^{(n,\alpha)}_\infty=\theta$, the $\{Z_i\}_{i\geq 1}$ are conditionally independent and identically distributed with distribution $\text{Bernoulli}(\theta)$, so that we have the integral representation
$$
  P\{Z_1=z_1,\dots,Z_m=z_m\} = \int_{[0,1]} \theta^{\sum_{i=1}^m z_i}(1-\theta)^{m-\sum_{i=1}^m z_i}\,d\mu(\theta)\,,
$$
for $(z_1,\dots,z_m)\in\{0,1\}^m$. For $k=0,1,\dots,m$, the event $\{\sum_{i=1}^m Z_i=k\}$ is the union of $\binom{m}{k}$ mutually exclusive and, by exchangeability, equiprobable events of the form $\{Z_1=z_1,\dots,Z_m=z_m\}$, in which $(z_1,\dots,z_m)\in\{0,1\}^m$, with $\sum_{i=1}^m z_i=k$. Therefore, it follows from the integral representation above and Definition \ref{dfn:coverage} that
\begin{equation*}
  P\left(C^{(n,\alpha)}_m=\frac{k}{m}\right) = P\!\left({\textstyle \sum_{i=1}^m Z_i=k}\right) = \binom{m}{k} \int_{[0,1]} \theta^k(1-\theta)^{m-k}\,d\mu(\theta)\,. \tag{$\dagger$}
\end{equation*}
Let the distribution $\mu$ of $C^{(n,\alpha)}_\infty$ be dominated by Lebesgue measure $\lambda$ with Radon-Nikodym derivative
\begin{equation*}
  (d\mu/d\lambda)(\theta) = \left(\frac{n!}{(b-1)!(g-1)!}\right)\theta^{b-1}(1-\theta)^{g-1}\,I_{[0,1]}(\theta)\,, 
\end{equation*}
up to almost everywhere $[\lambda]$ equivalence, in which $b=\lceil(1-\alpha)(n+1)\rceil$ and $g=n-b+1=\lfloor\alpha(n+1)\rfloor$. This is a version of the density of a random variable with $\mathrm{Beta}(b,g)$ distribution. Using $(\dagger)$ and the Leibniz rule for Radon-Nikodym derivatives (see \cite{schervish}, Theorem A.79), we have that
\begin{align*}
  P\left(C^{(n,\alpha)}_m=\frac{k}{m}\right) &= \binom{m}{k}\left(\frac{n!}{(b-1)!(g-1)!}\right)\int_{[0,1]}\theta^{k+b-1}(1-\theta)^{m-k+g-1}\,d\lambda(\theta) \\
  &= \binom{m}{k}\left(\frac{n!}{(b-1)!(g-1)!}\right)\left(\frac{(k+b-1)!(m-k+g-1)!}{(m+n)!}\right)\,.
\end{align*}
Since de Finetti's representation theorem states that the distribution $\mu$ of $C^{(n,\alpha)}_\infty$ is unique and that $(1/m)\sum_{i=1}^m Z_i$ converges almost surely to $C^{(n,\alpha)}_\infty$, when the batch size $m$ tends to infinity, the result follows by inspection of the distribution of the empirical coverage $C^{(n,\alpha)}_m$ in Theorem \ref{thm:finite}.
\end{proof}

\bibliographystyle{elsarticle-num} 

\bibliography{bibliography}

\end{document}